\documentclass[10pt,leqno,twoside]{amsart}
\usepackage{graphicx}

\newtheorem{theorem}{Theorem}[section]
\newtheorem{lemma}[theorem]{Lemma}
\newtheorem{proposition}[theorem]{Proposition}

\newtheorem{corollary}[theorem]{Corollary}
\theoremstyle{definition}
\newtheorem{definition}[theorem]{Definition}

\newtheorem{remark}[theorem]{Remark}

\numberwithin{equation}{section}
\usepackage{amsfonts}
\usepackage{multicol}

\usepackage{amsmath}
\usepackage{amssymb}

\begin{document}

%%%%%%%%%%%%%%%%%%%%%%%%%%%%%MACRO%%%%%%%%%%%%%%%%
\newcommand{\simge}{\ba{cc}\vspace*{-2.4mm}>\\ \sim\ea }
\newcommand{\simle}{\ba{cc}\vspace*{-2.4mm}<\\ \sim\ea }
\newcommand{\Cdot}{\!\cdot\!}
\newcommand{\sq}{{$\sqcap\!\!\!\!\sqcup$}}
\newcommand{\Eu}{{\rm I\,\!\! E}}
\newcommand{\Io}{\Int{\Omega}{}}
\newcommand{\Id}{\Int{\cald}{}}
\newcommand{\Div}{\mbox{\rm div}\,}
\newcommand{\tr}{\mbox{\rm tr}\,}
\newcommand{\grad}{\mbox{\rm grad}\,}
\newcommand{\supp}{\mbox{\rm supp}\,}
\newcommand{\curl}{\mbox{\rm curl}\,}
\newcommand{\Ido}{\Int{\partial\Omega}{}}
\newcommand{\IdS}{\Int{\Sigma}{}}
\newcommand{\Oint}[2]{{\displaystyle \oint_{#1}^{#2}}}
\newcommand{\Int}[2]{{\displaystyle \int_{ #1}^{ #2}}}
\newcommand{\Lim}[1]{{\displaystyle \lim_{ #1}}}
\newcommand{\Limsup}[1]{{\displaystyle \limsup_{\footnotesize #1}}}
\newcommand{\Liminf}[1]{{\displaystyle \liminf_{\footnotesize #1}}}
\newcommand{\Sup}[1]{{\displaystyle \sup_{#1}}}
\newcommand{\Inf}[1]{{\displaystyle \inf_{#1}}}
\newcommand{\Max}[1]{{\displaystyle \max_{#1}}}
\newcommand{\Min}[1]{{\displaystyle \min_{#1}}}
\newcommand{\Sum}[2]{{\displaystyle \sum_{#1}^{#2}}}
\newcommand{\Prod}[2]{{\displaystyle \prod_{#1}^{#2}}}
\newcommand{\BCup}[2]{{\displaystyle \bigcup_{#1}^{#2}}}
\newcommand{\BCap}[2]{{\displaystyle \bigcap_{#1}^{#2}}}
\newcommand{\Frac}[2]{\displaystyle{\frac{\displaystyle{#1}}{\displaystyle{#2}}}}
\newcommand{\norm}[1]{\left\|{#1}\right\|}
\newcommand{\Norm}[1]{\langle\langle{#1}\rangle\rangle_q}
\newcommand{\No}[1]{\langle\!\langle{#1}\rangle\!\rangle}
\newcommand{\NO}[1]{{\langle{#1}\rangle}_{\lambda,q}}
\newcommand{\beea}{\begin{eqnarray}}
\newcommand{\eeea}{\end{eqnarray}}
\newcommand{\ms}{\medskip\smallskip}
\newcommand{\bs}{\bigskip}
\newcommand{\ps}{\par\smallskip}
\newcommand{\bfe}{{\mbox{\boldmath $e$}} }
\newcommand{\pni}{{\par\noindent}}
\newcommand{\bfq}{{\mbox{\boldmath $q$}} }
\newcommand{\bfz}{{\mbox{\boldmath $z$}} }
\newcommand{\0}{{\mbox{\boldmath $0$}} }
\newcommand{\LE}{\!\!\!&\le&\!\!\!}
\newcommand{\BL}[1]{{\par\smallskip{\bf Lemma #1.}}}
\newcommand{\BT}[1]{{\par\smallskip{\bf Theorem #1.}}}
\newcommand{\Ln}{[\!|}
\newcommand{\Rn}{|\!]}
\newcommand{\n}[1]{{\Ln{#1}\Rn}} 
\newcommand{\nq}[1]{{\Ln{#1}\Rn}_{q}} 
\newcommand{\nqr}[1]{{\Ln{#1}\Rn}_{q,r}} 
\newcommand{\Nq}[1]{{\langle{#1}\rangle}_{q}} 
\newcommand{\Nql}[1]{{\langle{#1}\rangle}_{\lambda,q}} 
\newcommand{\Nqr}[1]{{\langle{#1}\rangle}_{q,r}}
\newcommand{\N}[1]{{|\!\!|\!\!|\,{#1}\,|\!|\!\!|_2}}
\newcommand{\BF}{\begin{footnotesize}}
\newcommand{\EF}{\end{footnotesize}}
\setlength{\jot}{.15in}
\newcommand{\pde}[2]{{\displaystyle \frac{\mbox{$\partial #1$}}{\mbox{$\partial #2$}}}}
\newcommand{\ode}[2]{{\displaystyle \frac{\mbox{$d #1$}}{\mbox{$d #2$}}}}
\newcommand{\f}[2]{\frac{\mbox{$#1$}}{\mbox{$ #2$}}}
\newcommand{\bi}{\begin{itemize}}
\newcommand{\ei}{\end{itemize}}
\newcommand{\ed}{\end{document}}
\newcommand{\be}{\begin{equation}}
\newcommand{\ba}{\begin{array}}
\newcommand{\ea}{\end{array}}
\newcommand{\ee}{\end{equation}}
\newcommand{\eeq}[1]{\label{eq:#1}\end{equation}}

\newcommand{\R}{{\mathbb R}}
\newcommand{\bH}{{\mathbb H}}

\newcommand{\compl}{{\mathbb C}}
\def\Id{\mbox{\boldmath $1$}}
\def\zero{\mbox{\boldmath $0$}}
\newcommand{\PP}{{\rm I\!\!\,P}}
\newcommand{\nat}{{\mathbb N}}
\newcommand{\bfpsi}{\mbox{\boldmath $\psi$}}
\newcommand{\bfchi}{\mbox{\boldmath $\chi$}}
\newcommand{\bfomega}{\mbox{\boldmath $\omega$}}
\newcommand{\bfvaromega}{\mbox{\boldmath $\varpi$}}
\newcommand{\bfOmega}{\mbox{\boldmath $\Omega$}}
\newcommand{\bfTheta}{\mbox{\boldmath $\Theta$}}
\newcommand{\bfxi}{\mbox{\boldmath $\xi$}}
\newcommand{\bfmu}{\mbox{\boldmath $\mu$}}
\newcommand{\bfx}{\mbox{\boldmath $x$}}
\newcommand{\bfy}{\mbox{\boldmath $y$}}
\newcommand{\bfPsi}{\mbox{\boldmath $\Psi$}}
\newcommand{\bfphi}{\mbox{\boldmath $\varphi$}}
\newcommand{\bfhi}{\mbox{\boldmath $\phi$}}
\newcommand{\bfPhi}{\mbox{\boldmath $\Phi$}}
\newcommand{\bfv}{{\mbox{\boldmath $v$}} }
\newcommand{\bfu}{{\mbox{\boldmath $u$}} }
\newcommand{\bfsf}{{\mbox{\footnotesize\boldmath $s$}} }
\newcommand{\bfuf}{{\mbox{\footnotesize\boldmath $u$}} }
\newcommand{\bfw}{{\mbox{\boldmath $w$}} }
\newcommand{\bff}{{\mbox{\boldmath $f$}} }
\newcommand{\bfa}{{\mbox{\boldmath $a$}} }
\newcommand{\bfi}{{\mbox{\boldmath $i$}} }
\newcommand{\bfj}{{\mbox{\boldmath $j$}} }
\newcommand{\bfc}{{\mbox{\boldmath $c$}} }
\newcommand{\bfo}{{\mbox{\boldmath $o$}} }
\newcommand{\bfp}{{\mbox{\boldmath $p$}} }
\newcommand{\bfkp}{{\mbox{\footnotesize{\boldmath $k$}}} }
\newcommand{\bfka}{{\mbox{\footnotesize{\boldmath $k^*$}}} }
\newcommand{\bft}{{\mbox{\boldmath $t$}} }
\newcommand{\bfd}{{\mbox{\boldmath $d$}} }
\newcommand{\bfl}{{\mbox{\boldmath $l$}} }
\newcommand{\bfr}{{\mbox{\boldmath $r$}} }
\newcommand{\bfk}{{\mbox{\boldmath $k$}} }
\newcommand{\bfA}{{\mbox{\boldmath $A$}} }
\newcommand{\bfS}{{\mbox{\boldmath $S$}} }
\newcommand{\bfO}{{\mbox{\boldmath $O$}} }
\newcommand{\bfM}{{\mbox{\boldmath $M$}} }
\newcommand{\bfP}{{\mbox{\boldmath $P$}} }
\newcommand{\bfB}{{\mbox{\boldmath $B$}} }
\newcommand{\bfR}{{\mbox{\boldmath $R$}} }
\newcommand{\bfC}{{\mbox{\boldmath $C$}} }
\newcommand{\bfD}{{\mbox{\boldmath $D$}} }
\newcommand{\bfQ}{{\mbox{\boldmath $Q$}} }
\newcommand{\bfZ}{{\mbox{\boldmath $Z$}} }
\newcommand{\bfG}{{\mbox{\boldmath $G$}} }
\newcommand{\bfE}{{\mbox{\boldmath $E$}} }
\newcommand{\bfX}{{\mbox{\boldmath $X$}} }
\newcommand{\bfY}{{\mbox{\boldmath $Y$}} }
\newcommand{\bfH}{{\mbox{\boldmath $H$}} }
\newcommand{\bfI}{{\mbox{\boldmath $I$}} }
\newcommand{\bfJ}{{\mbox{\boldmath $J$}} }
\newcommand{\bfN}{{\mbox{\boldmath $N$}} }
\newcommand{\bfh}{{\mbox{\boldmath $h$}} }
\newcommand{\bfm}{{\mbox{\boldmath $m$}} }
\newcommand{\bfone}{{\mbox{\boldmath $1$}} }
\newcommand{\hs}{{\rm I}\!\!\,{\rm R}^3_+}
\newcommand{\cala}{{\cal A}}
\newcommand{\calb}{{\cal B}}
\newcommand{\calc}{{\cal C}}
\newcommand{\cald}{{\cal D}}
\newcommand{\cale}{{\cal E}}
\newcommand{\calf}{{\cal F}}
\newcommand{\calg}{{\cal G}}
\newcommand{\calh}{{\cal H}}
\newcommand{\cali}{{\cal I}}
\newcommand{\calj}{{\cal J}}
\newcommand{\calk}{{\cal K}}
\newcommand{\call}{{\cal L}}
\newcommand{\calm}{{\cal M}}
\newcommand{\caln}{{\cal N}}
\newcommand{\calo}{{\cal O}}
\newcommand{\calp}{{\cal P}}
\newcommand{\calq}{{\cal Q}}
\newcommand{\calr}{{\cal R}}
\newcommand{\cals}{{\cal S}}
\newcommand{\cT}{{\mathcal T}}
\newcommand{\calu}{{\cal U}}
\newcommand{\calv}{{\cal V}}
\newcommand{\calx}{{\cal X}}
\newcommand{\caly}{{\cal Y}}
\newcommand{\calw}{{\cal W}}
\newcommand{\calz}{{\cal Z}}
\newcommand{\bfsigma}{\mbox{\boldmath $\sigma$}}
\newcommand{\bfSigma}{\mbox{\boldmath $\Sigma$}}
\newcommand{\bftau}{\mbox{\boldmath $\tau$}}
\newcommand{\bfeta}{\mbox{\boldmath $\eta$}}
\newcommand{\bfT}{{\mbox{\boldmath $T$}} }
\newcommand{\bfV}{{\mbox{\boldmath $V$}} }
\newcommand{\bfU}{{\mbox{\boldmath $U$}} }
\newcommand{\bfW}{{\mbox{\boldmath $W$}} }
\newcommand{\bfF}{{\mbox{\boldmath $F$}} }
\newcommand{\bfK}{{\mbox{\boldmath $K$}} }
\newcommand{\bfL}{{\mbox{\boldmath $L$}} }
\newcommand{\bfb}{{\mbox{\boldmath $b$}} }
\newcommand{\bfg}{{\mbox{\boldmath $g$}} }
\newcommand{\bfn}{{\mbox{\boldmath $n$}} }
\newcommand{\bfs}{{\mbox{\boldmath $s$}} }
\newcommand{\cf}{{\it cf.} }
\newcommand{\io}{\int_\Omega}
\newcommand{\1}{\item[({\it i})]}
\newcommand{\2}{\item[({\it ii})]}
\newcommand{\3}{\item[({\it iii})]}
\newcommand{\4}{\item[({\it iv})]}
\newcommand{\5}{\item[({\it v})]}
\newcommand{\6}{\item[({\it vi})]}
\newcommand{\7}{\item[({\it vii})]}
\newcommand{\8}{\item[({\it viii})]}
\newcommand{\9}{\item[({\it xi})]}
\newcommand{\ido}{\int_{\partial\Omega}}
\newcommand{\half}{\mbox{$\frac{1}{2}$}}
\def\parallel{\|}
\def\mid{|}
\def\Bbb R{\real}
\def\hat{\widehat}
\def\tilde{\widetilde}
\def\bar{\overline}
\newcommand{\threehalves}{3\over 2}
\newcommand{\bfPi}{\mbox{\boldmath $\Pi$}}
\newcommand{\bfXi}{\mbox{\boldmath $\Xi$}}
\newcommand{\bfalpha}{\mbox{\boldmath $\alpha$}}
\newcommand{\bfbeta}{\mbox{\boldmath $\beta$}}
\newcommand{\bfgamma}{\mbox{\boldmath $\gamma$}}
\newcommand{\bfdelta}{\mbox{\boldmath $\delta$}}
\newcommand{\bfzeta}{\mbox{\boldmath $\zeta$}}
\newcommand{\bfUpsilon}{\mbox{\boldmath $\Upsilon$}}
\newcommand{\bfGamma}{\mbox{\boldmath $\Gamma$}}
\newcommand{\bfcala}{\mbox{\boldmath ${\cal A}$}}
\newcommand{\bfcalm}{\mbox{\boldmath ${\cal M}$}}
\newcommand{\bfcaln}{\mbox{\boldmath ${\cal N}$}}
\newcommand{\bfcalq}{\mbox{\boldmath ${\cal Q}$}}
\newcommand{\bfcalb}{\mbox{\boldmath ${\cal B}$}}
\newcommand{\bfcalc}{\mbox{\boldmath ${\cal C}$}}
\newcommand{\bfcali}{\mbox{\boldmath ${\cal I}$}}
\newcommand{\bfcalg}{\mbox{\boldmath ${\cal G}$}}
\newcommand{\bfcalh}{\mbox{\boldmath ${\cal H}$}}
\newcommand{\bfcalk}{\mbox{\boldmath ${\cal K}$}}
\newcommand{\bfcalt}{\mbox{\boldmath ${\cal T}$}}
\newcommand{\bfcalx}{\mbox{\boldmath ${\cal X}$}}
\newcommand{\bfcall}{\mbox{\boldmath ${\cal L}$}}
\newcommand{\bfcalf}{\mbox{\boldmath ${\cal F}$}}
\newcommand{\bfcalr}{\mbox{\boldmath ${\cal R}$}}
\newcommand{\bfcals}{\mbox{\boldmath ${\cal S}$}}
\newcommand{\bfcalw}{\mbox{\boldmath ${\cal W}$}}
\newcommand{\bfcalu}{\mbox{\boldmath ${\cal U}$}}
\newcommand{\bfcalv}{\mbox{\boldmath ${\cal V}$}}
\newcommand{\bfcalz}{\mbox{\boldmath ${\cal Z}$}}
\pagenumbering{roman}
\newcommand{\art}[6]{{\I[{\sc #1,}] {#2}, {\it #3}, {\bf #4}, {#5} {[#6]}}}
\newcommand{\ED}{\end{description}}
\newcommand{\I}{\item }
\newcommand{\ra}{\rm a}
\newcommand{\rb}{\rm b}
\newcommand{\rc}{\rm c}
\newcommand{\Hsp}{{\rm I}\!\!\,{\rm R}^n_+}
\newcommand{\Hsn}{{\rm I}\!\!\,{\rm R}^n_-}
\newcommand{\po}[1]{\mbox{$\displaystyle \frac{\mbox{$\partial #1$}}
{\mbox{$\partial x_{1}$}}$}}
\newcommand{\PO}[1]{\mbox{$\displaystyle \frac{\mbox{$\partial #1$}}
{\mbox{$\partial y_{1}$}}$}}
\newcommand{\OP}{\left(\Delta+2\lambda\PO{}\right)}
\newcommand{\op}{\left(\Delta+2\lambda\po{}\right)}
\newcommand{\ft}[1]{
\Frac{1}{(2\pi)^{n/2}}\Int{{\Bbb R}^{n}}{}e^{i{\bf x}\cdot \bfxi}
#1(\xi)d\xi}
\newcommand{\Ft}[1]{
\Frac{1}{2\pi}\Int{{\Bbb R}^{2}}{}e^{i{x}\cdot \xi}
#1(\xi)d\xi}
\newcommand{\Z}{\item[({\it a})]}
\newcommand{\B}{\item[({\it b})]}
\newcommand{\C}{\item[({\it c})]}
\newcommand{\D}{\item[({\it d})]}
\newcommand{\E}{\item[({\it e})]}
\newcommand{\Footnote}{~\footnote}
\newcommand{\ie}{{\it i.e.}}
\newcommand{\dist}{\mbox{\rm dist\,}}
\newcommand{\const}{\mbox{\rm const}}
\newcommand{\trace}{\mbox{\rm trace}}
\newcommand{\Bo}{\par\hfill{$\Box$}\par\noindent}
\newcommand{\Nor}[1]{\langle{#1}\rangle_q}
\newcommand{\vs}{\vspace*{.5cm}\par\noindent}
\newcommand{\Vs}{\vspace*{.6cm}\par\noindent}
\newcommand{\Vvs}{\vspace*{.7cm}\par\noindent}
\newcommand{\VVs}{\vspace*{.8cm}\par\noindent}
\newcommand{\Br}{\begin{remark}\begin{rm}}
\newcommand{\Er}{\end{rm}\end{remark}}
\newcommand{\Bp}{\begin{proposition}\begin{sl}}
\newcommand{\EP}[1]{\end{sl}\label{proposition:#1}\end{proposition}}
\newcommand{\Bt}{\begin{theorem}\begin{sl}}
\newcommand{\Et}{\end{sl}\end{theorem}}
\newcommand{\Bl}{\begin{lemma}\begin{sl}}
\newcommand{\Bc}{\begin{corollary}\begin{sl}}
\newcommand{\Ec}{\end{sl}\end{corollary}}
\newcommand{\ET}[1]{\end{sl}\label{theorem:#1}\end{theorem}}
\newcommand{\EDD}[1]{\end{rm}\label{definition:#1}\end{definition}}
\newcommand{\EL}[1]{\end{sl}\label{lemma:#1}\end{lemma}}
\newcommand{\ER}[1]{\end{rm}\label{remark:#1}\end{remark}}
\newcommand{\EC}[1]{\end{sl}\label{corollary:#1}\end{corollary}}
\newcommand{\essup}[1]{{\rm ess}\,{{\displaystyle \sup_{\hspace*{-5mm}{#1}}}}}
%%%%%%%%%%%%%%%%%%%%%%%ENDOFMACRO%%%%%%%%%%%%%%%%%

\pagenumbering{arabic}
\newcommand{\QED}{{\par\hfill$\square$\par}}
\newcommand{\NP}{{\nabla_H \,}}
\newcommand{\DP}{{\Div_H \,}}

\title{Strong Time-Periodic Solutions to the 3D Primitive Equations subject to  Arbitrary Large Forces}

\author[Giovanni P. Galdi]{Giovanni P. Galdi}
\address{Benedum Engineering Hall, University of Pittsburgh, Pittsburgh, PA 15261, USA} 
\email{galdi@pitt.edu}

\author[Matthias Hieber]{Matthias Hieber}
\address{Department of Mathematics, TU Darmstadt, Schlossgartenstr. 7, 64289 Darmstadt, Germany, and 
Benedum Engineering Hall, University of Pittsburgh, Pittsburgh, PA 15261, USA}  
\email{hieber@mathematik.tu-darmstadt.de}

\author[Takahito Kashiwabara]{Takahito Kashiwabara}
\address{Department of Mathematics, Tokyo Institute of Technology, 2-12-1 Ookayama, Meguro, Tokyo 152-8511, Japan}
\email{tkashiwa@math.titech.ac.jp}


\thanks{*The work of G.P.Galdi was partially supported by  NSF DMS Grant-1311983, and a Mercator Professorship at the Department of Mathematics TU Darmstadt.}

\begin{abstract}
We show that the three-dimensional primitive equations admit a strong time-periodic solution of period  $T>0$, provided the forcing term 
$f\in L^2(0,\cT; L^2(\Omega))$ is a   time-periodic function of the same period. No restriction on the magnitude of $f$ is assumed. As a corollary, if, in particular, $f$ is 
time-independent, the corresponding solution is steady-state.  
\end{abstract}

\maketitle

\section{Introduction}
Consider the primitive equations in the isothermal setting, i.e.  assuming that the temperature $\theta$ equals a constant $\theta_0$. In this case, the primitive equations consist of the 
following set of equations   
\begin{equation}\label{eq:1.1}
\begin{array}{rll}
\partial_t v + u\cdot\nabla v - \Delta v + \nabla_H\pi & =  f \quad   & \text{ in } \Omega\times(0,\cT), \\
    		\partial_z\pi & =  0  \quad &\text{ in } \Omega\times(0,\cT), \\
		\mathrm{div}\,u & =  0  \quad &\text{ in } \Omega\times(0,\cT), \\
                             v(0) & = a. & 
\end{array}
\end{equation}
Here $\Omega = G \times(-h,0)$, where $G=(0,1)^2$, $h>0$, and $\cT>0$ . The velocity $u$ of the fluid is given by $u=(v,w)$ with $v=(v_1,v_2)$, and where $v$ and $w$ denote  
the  horizontal and vertical components of $u$, respectively. Furthermore, $\pi$ denotes  the pressure of the fluid (more precisely, $\pi = p + \theta_0z$, where $p$ is the original 
pressure, $z\in(-h,0)$) and $f$ a given external force. 
The symbol $\nabla_H=(\partial_x,\partial_y)^\top$  denotes the horizontal gradient, $\Delta$ the  three-dimensional Laplacian and 
$\nabla$ and $\mathrm{div}$ the three dimensional gradient and divergence operators. The above equations take into account, by scale analysis,  the hydrostatic approximation of the 
Navier-Stokes equations; for more details see e.g. \cite{TZ04}, \cite{Val06}.  

The system is complemented by the set of boundary conditions
\begin{equation}\label{eq:1.2}
\begin{array}{rll}
\partial_z v & =   0, \quad w=0  &\text{ on }\; \Gamma_u\times(0,\cT),  \\
	   v & = 0, \quad w=0 & \text{ on }\; \Gamma_b\times(0,\cT), \\
&	\text{$u$, $\pi$ are periodic}   & \text{ on }\; \Gamma_l\times(0,\cT).  \\
\end{array}
\end{equation}
Here   $\Gamma_u := G\times\{0\}$, $\Gamma_b := G \times\{-h\}$, $\Gamma_l :=\partial G\times[-h,0]$ 
denote the upper, bottom and  lateral parts of the boundary $\partial\Omega$, respectively.

The full primitive equations were introduced and investigated for the first time by Lions, Temam and Wang in \cite{LTW92, LTW95}.  They proved the 
existence of a global weak solution for this set of equations for initial data $a \in L^2$. The existence of a local, strong solution with data $a \in H^1$ was proved first 
by Guill\'en-Gonz\'alez, Masmoudi and Rodiguez-Bellido in \cite{GMR01}. 

In 2007, Cao and Titi \cite{CT07} proved a breakthrough result for this set of equations which states, roughly speaking,  that there exists a unique, {\em global strong} solution to 
the primitive equations for {\it arbitrary} initial data $a\in H^1$. Note that the boundary conditions on $\Gamma_b\cup\Gamma_l$ considered there are different from the ones 
we are imposing in \eqref{eq:1.2}. Successively, in \cite{KZ07} Kukavica and Ziane considered  the primitive equations subject to boundary conditions 
as in \eqref{eq:1.2}, and proved global strong well-posedness of the primitive equations with respect to arbitrary $H^1$-data. For different approaches see also 
Kobelkov \cite{Kob06} and Kukavica, Pei, Rusin and Ziane \cite{KPRZ14}. 
 
It is worth emphasizing that, while the fundamental problem of global existence and uniqueness for the initial-value problem can be considered to a great extent settled, at least in the $L^2$ framework, other important issues like existence of {\em strong} time-periodic (and, in particular, steady-state) solutions to the primitive equations appear to be at a stage where  further investigation is still required.     
In this regard, we recall that the question of whether system \eqref{eq:1.1} admits time-periodic solutions  was first addressed by Tachim Medjo \cite{Tac10}.\footnote{Even in the more general non-isothermal framework.} There, under the assumption of {\em ``small'' forcing term}, existence (and uniqueness) of strong solutions is achieved by the classical Galerkin method. More recently, 
Hsia and Shiue \cite{HS13}  proved a similar result by a different method suggested by Serrin \cite{Ser59}, again under a suitable {\em smallness condition} on the forcing term. Furthermore, they showed asymptotic stability of such solutions when the initial perturbations are sufficiently small.  Notice that, as  corollary, both results in \cite{Tac10,HS13} furnish existence of steady-state solutions to the primitive equations for forcing terms of {\em restricted} magnitude.  

At this point, it must be observed that the smallness assumption on the forcing term is undesired and, most of all, appears somehow unexpected if one compares this situation with the classical Navier-Stokes theory. There, even though the initial-value problem still lacks of a global existence result for strong solutions with initial data of arbitrary size, nevertheless the steady-state boundary-value problem is known to have a smooth solution for (smooth) forcing term of arbitrary 
magnitude since the fundamental work of J.Leray. However, it should be  added, also in the light of  the contributions \cite{Tac10,HS13}, that the achievement of a result of this type for the primitive equations does not seem to be obvious, at least if one uses the classical methods employed for the Navier-Stokes equations.       

The main objective of this article is to prove  existence of {\em strong} time-periodic  solutions to the primitive equations of the form \eqref{eq:1.1}$_{1,2,3}$ for {\em arbitrary} (time-periodic) $f \in L^2(0,\cT,L^2(\Omega))$, hereby {\em without 
assuming any smallness condition on $f$}. As a byproduct, this result provides an analogous one for steady-state solutions. 

As we hinted earlier on, the approach  we use differs from ``standard'' ones, and   is based on the 
following three steps: First, we construct a (suitable) {\em weak} time-periodic solution, $v$, to \eqref{eq:1.1}$_{1,2,3}$ corresponding to the given $f$, 
by combining  classical Galerkin's method with Brouwer's fixed point theorem. Secondly, we show the existence of a unique, {\em strong} solution $u$ to the {\em initial-value} problem  \eqref{eq:1.1} for arbitrary $f\in L^2(0;\cT; L^2(\Omega))$, and $a$ in a subspace of $H^1(\Omega)$, by using the arguments of \cite{HK15}. Finally, we look at $v$ as a weak solution to the {\em initial-value} problem and  employ  a weak-strong uniqueness result of the type proved by  Guill\'en-Gonz\'alez,  Masmoudi  and Rodriguez-Bellido 
 in \cite{GMR01}, which then implies $v\equiv u$, thus furnishing the main 
results of this article stated as Theorem \ref{thm:main} and Corollary \ref{cor}.  

The plan of the paper is the following. In Section 2 we make some preliminary considerations and give the statement of our main results (Theorem \ref{thm:main} and Corollary \ref{cor}). In Section 3 we then show the existence of a weak time-periodic solution corresponding to forcing terms of arbitrary size. The following Section 4 is dedicated to the proof of existence and uniqueness of (an equivalent form of) the initial-value problem \eqref{eq:1.1} for arbitrary $f$ and $a$ in appropriate function classes. Finally, in Section 5, we give the proof of Theorem \ref{thm:main}.

\section{Preliminaries and Main Results}
Following Lions, Temam and Wang \cite{LTW92, LTW95} and Cao and Titi \cite{CT07}, we rewrite  the primitive equations given in \eqref{eq:1.1} subject 
to the boundary conditions \eqref{eq:1.2}  in the following equivalent form. Since the vertical component $w$ of $u$ is determined by the incompressibility condition we have  
\begin{equation*} %\label{eqp.1}
	w(x,y,z) = \int_{z}^0 \mathrm{div}_H\,v(x,y,\zeta)\,d\zeta, \qquad (x,y)\in G,\; -h<z<0,
\end{equation*}
due to the boundary condition  $w=0$ on $\Gamma_u$. The further boundary condition $w=0$ on $\Gamma_b$ gives rise to the constraint
\begin{equation*}
	\mathrm{div}_H\,\bar v = 0 \quad\text{in}\quad G,
\end{equation*}
where $\bar v$ stands for the average of $v$ in the vertical direction, i.e.,
\begin{equation} \label{eqp.1}
	\bar v(x,y) := \frac1h\int_{-h}^0 v(x,y,z)\,dz, \quad (x,y) \in G.
\end{equation}
Then problem \eqref{eq:1.1}-\eqref{eq:1.2} is equivalent to finding a function $v:\Omega\to\mathbb R^2$ and a function $\pi:G\to\mathbb R$ satisfying the set of equations 
\begin{equation} \label{eqp.2}
	\begin{array}{rll}
		\partial_tv + v\cdot\nabla_H v + w(v)\,\partial_zv - \Delta v + \nabla_H\pi  & = f & \text{ in } \Omega\times(0,T), \\
		w(v) & = \int_z^0 \mathrm{div}_H\,v\,d\zeta & \text{ in } \Omega\times(0,T), \\
		\mathrm{div}_H\,\bar v &= 0 & \text{ in } G\times(0,T), \\
                                    v(0) & = a,&  
	\end{array}
\end{equation}
as well as  the boundary conditions
\begin{equation} \label{eqp.3}
	\begin{array}{rll}
		\partial_{z} v & = 0  & \text{ on } \Gamma_u\times(0,T), \\
		v &= 0 & \text{ on } \Gamma_b\times(0,T), \\
	&	\text{$v$ and $\pi$ are periodic}  & \text{ on } \Gamma_l\times(0,T). \\
	\end{array}
\end{equation}
The following terminology for  describing the periodic boundary conditions will be useful. Let $m \in \{0,1\}$. We then say that a smooth function 
$f: \overline {\Omega} \to \R$ is {\em space periodic of order $m$ on $\Gamma_l$} if
\begin{equation*}
	\frac{\partial^\alpha f}{\partial x^\alpha}(0,y,z) = 
\frac{\partial^\alpha f}{\partial x^\alpha}(1,y,z) \, \mbox{ and } \, \frac{\partial^\alpha f}{\partial y^\alpha}(x,0,z) = \frac{\partial^\alpha f}{\partial y^\alpha}(x,1,z),
\end{equation*}
for all $\alpha=0,\dots,m$. Note that we do not consider any symmetry conditions in the $z$-direction.
The Sobolev spaces equipped with space-periodic boundary conditions in the horizontal directions are defined by\footnote{Throughout the paper we shall use the same font style to denote scalar, vector and
tensor--valued functions and corresponding function spaces.}
\begin{align*}
	H^{m}_{\mathrm{per}}(\Omega) &:= \{f\in H^{m}(\Omega):  \text{$f$ is space-periodic of order $m-1$ on $\Gamma_l$} \}. 
\end{align*}
Note that $C^\infty_{\mathrm{per}}(\overline\Omega) := \{f\in C^\infty(\overline\Omega): \text{$f$ is space-periodic of arbitrary order on $\Gamma_l$} \} $ is dense in 
$H^{m}_{\mathrm{per}}(\Omega)$. We now introduce the function spaces $\bH, \bH^1$ and $\bH^2$ by  
\begin{align*}
\bH(\Omega)    &:=\{v \in L^2_{\mathrm{per}}(\Omega): \DP \bar v = 0 \} \\
\bH^1(\Omega)  &:=\{v \in H^1_{\mathrm{per}}(\Omega): \DP \bar v =0, v=0 \mbox{ on } \Gamma_b \} \\    
\bH^2(\Omega)  &:=\{v\in \bH^1(\Omega)\cap H^2(\Omega)_\text{{per}}: {\partial_z v}=0 \mbox{ on } \Gamma_u \}.
\end{align*}
and denote its norms by  $\|\cdot\|_2$, {$\|\cdot\|_{\bH^1}$, $\|\cdot\|_{\bH^2}$}, respectively. As usual, $(\cdot,\cdot)$ stands for the usual $L^2$ scalar product. 
 
Moreover, given an interval $I \subset \R$, the space $C_w(I;\bH(\Omega))$ stands for the class of functions $v:I \to \bH(\Omega)$ such that $t \mapsto (v(t),\psi)$ is 
continuous for all $\psi\in \bH(\Omega)$. 
\par
Finally, we say that a function $f\in L^1(0,\cT;L^1(\Omega))$, all $\cT>0$,  is $T$-periodic ($T>0$) if $f(t,x)=f(t+T,x)$ for a.a. $(t,x)\in [0,\infty)\times \Omega$.
\begin{definition}\label{def:weak} Let $T>0$ and $f\in L^1(0,\cT;L^2(\Omega))$,  all $\cT>0$. A function  $v:[0,\infty)\times \Omega \to \R^2$ is called  
{\em weak $T$-periodic solution} to \eqref{eqp.2}$_{1,2,3}$  and \eqref{eqp.3} if 
\begin{itemize}
\item[i)] $v\in C_w([0,\cT];\bH(\Omega))\cap L^2(0,\cT;\bH^1(\Omega))$ for all $\cT >0$,  
\item[ii)] For all $\cT>0$ and all $\varphi\in C^1([0,\cT];\bH^1(\Omega))\cap L^2(0,\cT;\bH^2(\Omega))$
	
\begin{equation*}%\label{eq:2.1}
\Int0t\big\{(v,\partial_t\varphi)  - (\nabla v,\nabla\varphi)+(v\cdot\NP\varphi,v+w(v)\,{\partial_z}\varphi,v)\big\}
=-\Int0t(f,\varphi)+(v(t),\varphi(t))-(v(0),\varphi(0)), \quad t\in (0,\cT), 
\end{equation*}
\item[iii)] For all $\cT>0$, all $t\in (0,\cT]$ and a.a. $s\in[0,t)$, $v$ satisfies the {\em strong energy inequality}
$$  
\|v(t)\|_2^2+2 \int_s^t\|\nabla v(\tau)\|_2^2 d\tau\le \|v(s)\|_2^2+2\int_s^t(f(\tau),v(\tau)) d\tau, 
$$
\item[iv)] $v(t+T,x)=v(T,x)$ for all $t\ge 0$ and a.a. $x\in\Omega$.
\end{itemize}
A weak $T$-periodic solution $v$ is called a {\em strong} if, in addition to i)--iv), it holds  
$v\in C([0,\cT];\bH^1(\Omega))\cap L^2(0,\cT;\bH^2(\Omega))$, $\partial_tv\in L^1(0,\cT; L^2(\Omega))$, for all  $\cT>0$.
\par
If, in particular, $f\in \mathbb H(\Omega)$ is independent of $t\ge 0$, we say that $v_s\in \mathbb H^1(\Omega)$ is a {\em weak steady-state solution} to  \eqref{eqp.2}$_{1,2,3}$--\eqref{eqp.3}  if $v_s$ satisfies condition ii) for all $\varphi\in \mathbb H^2(\Omega)$. A weak steady-state solution is called {\em strong } if $v_s\in \mathbb H^2(\Omega)$.

\end{definition}

\begin{remark}
It is worth noticing that every term in ii) is well defined. This is obvious for all linear terms in $v$. As for the nonlinear ones, by H\"older's inequality and Sobolev embeddings we obtain for all $v_1,v_2\in H^1(\Omega)$ and all $v_3 \in H^2(\Omega)$ 
\begin{equation}\label{eq:2.2}
\Int0t|(v_1\cdot\NP v_3,v_2)|\le \Int0t\|v_1\|_3\|\NP v_3\|_2\|v_2\|_6\le C \int_0^t\|v_1\|_{2}^{\frac12} \|v_1\|_{H^1}^{\frac12} \|v_2\|_{H^1} \|\NP v_3\|_2\,.
\end{equation}
Furthermore, using inequality (93) from \cite{CT07}, we show that
\begin{equation}\label{eq:2.3}
\Int0t|(w(v_1)\,{\partial_z}v_3,v_2)|\le C\Int0t \|\NP v_1\|_{2}\|v_2\|_2^{\frac12}\| v_2\|_{H^1}^{\frac12} \|\partial_z v_3\|_2^{\frac12} \|\partial_z v_3\|_{H^1}^{\frac12} \,,
\end{equation}
which proves the claim.
We note that the proof of the latter does not use boundary conditions for $v_i$, $i=1,2,3$,  which in \cite{CT07} are different than those adopted here.

\end{remark}
\begin{remark} The above  definition 
of a weak time-periodic solution is somewhat different than the one typically found  in the literature. However, this formulation is needed when  we will compare these solutions with solutions to the initial-value problem; see Section 4.
\end{remark}
We are now in the position to state the main result of this article.

\begin{theorem}\label{thm:main}
Let $T>0$ and let  $f\in L^2(0,\cT;L^2(\Omega))$, all $\cT>0$, be $T$-periodic.   
Then  problem \eqref{eqp.2}$_{1,2,3}$--\eqref{eqp.3} has at least one  corresponding 
strong $T$-periodic solution. 
%$u_p$ satisfying  $u_p\in C([0,\cT]; \bH^1(\Omega)) \cap L^2(0,\cT; \bH^2(\Omega))$ for all $\cT>0$.  
\end{theorem}

The above result at once implies the following one.
\begin{corollary}\label{cor}
Let $f \in L^2(\Omega)$.   Then problem  \eqref{eqp.2}$_{1,2,3}$--\eqref{eqp.3} has at least one  
corresponding strong steady-state solution.
%there exists at least one  strong steady-state solution $u_{S} \in \bH^1(\Omega) \cap \bH^2(\Omega)$ to problem \eqref{eqp.2}--\eqref{eqp.3}.  
\end{corollary}

We emphasize at this point that, in contrast to previous known results \cite{HS13,Tac10}, our findings {\em do not require any smallness condition on $f$}.

\section{Weak Time-Periodic Solutions}
Objective of this section is to show that the class of weak $T$-periodic  solutions to \eqref{eqp.2}$_{1,2,3}$-- \eqref{eqp.3} is not empty under suitable assumptions on $f$. Precisely, we have the following.

\begin{proposition}\label{prop:3.1}
%Let $f:(0,\infty)\times \Omega \to \R^2$ be a time-periodic function of period $T>0$ such that  $f\in L^2(0,T;\bH(\Omega))$. 
Let $T>0$, and let  $f\in L^2(0,\cT;L^2(\Omega))$,  all $\cT>0$,  $T$-periodic. 
Then, there exists at least one weak $T$-periodic 
solution to \eqref{eqp.2}$_{1,2,3}$--\eqref{eqp.3}
\end{proposition}

\noindent
{\em Proof}. Even though our  definition of a weak time-periodic solution is somehow different than the one usually given in the  literature, the proof of its existence,  
based on  the Faedo-Galerkin method, is quite standard; see, e.g., \cite{Pro}, \cite[Chapter 4, \S 6.2]{Lio69}, \cite[pp. 256-260]{GS06}, \cite{Gal00}. 
For this reason, we shall  only give the main arguments, referring the reader to the above papers for further details. 

Let $\{\psi_n\}\subset \bH^2(\Omega)$ be an orthonormal basis of $\bH(\Omega)$ dense in $\bH^1(\Omega)$ and $\bH^2(\Omega)$. For example,  we may take the  eigenvectors of the 
hydrostatic Stokes operator in the $L^2$ setting (see \cite[\S 4]{HK15}). Let
$$
v_m:=\sum_{k=1}^mc_{mk}(t)\psi_k(x)
$$
where, for all $r=1,\ldots,m$,
\begin{equation}\label{eq:2.4}
\ode{}{t}(v_m,\psi_r)=(v_m\cdot\NP \psi_r,v_m)+(w(v_m)\,{\partial_z}\psi_r,v_m)+(\nabla v_m,\nabla\psi_r)+(f,\psi_r)\,. 
\end{equation}
We now show a uniform bound in time on the functions $c_{mk}(t)$ for  $k=1,\ldots,m$, which implies that the system \eqref{eq:2.4} has a solution 
$c(t):=(c_{m1}(t),\ldots,c_{mm}(t))$ for all $t\ge 0$ and all $m\ge 1$. In fact, multiplying both sides of \eqref{eq:2.4} by $c_{mr}(t)$, summing over $r$, and taking into account 
that by \eqref{eqp.2}$_{2,3,4}$
$$
(v_m\cdot\NP v_m,v_m)+(w(v_m)\,{\partial_z}v_m,v_m)=0\,,
$$
we get
\begin{equation}\label{eq:2.5}
\ode{}{t}\|v_m\|_2^2+2\|\nabla v_m(t)\|_2^2=2(f,v_m)\,.
\end{equation}
Since $v_m=0$ at $x_3=-h$, we have 
\be\label{eq:2.6}
\|v_m\|_2^2\le h^2\|{\partial_z}v_m\|_2^2, 
\ee
so that by Schwartz's inequality and \eqref{eq:2.5} we infer
$$
\ode{}{t}\|v_m\|_2+\frac2{h^2}\|v_m\|_2\le 2\|f\|_2\,.
$$
Integrating the latter from $t=0$ to arbitrary $t>0$ we get
\begin{equation}\label{eq:2.7}
{\rm e}^{2t/h^2}\|v_m(t)\|_2\le \|v_m(0)\|_2+ 2\int_0^t {\rm e}^{2\tau/h^2}\|f(\tau)\|_2d\tau\,,
\end{equation}
thus deducing the claimed uniform bound for $c(t)$, once we observe that, by the orthonormality property of $\{\psi_n\}$, $|c(t)|=\|v_m(t)\|_2$.

Next, choose $R>0$ such that  $R({\rm e}^{2t/h^2}-1) \ge 2\int_0^T {\rm e}^{2\tau/h^2}\|f(\tau)\|_2\,d\tau$ 
and let $\mathbb B_R^m$ the ball in  $\R^m$ centered at the origin with radius $R$. It follows from \eqref{eq:2.7} that $|c(T)|=\|v_m(T)\|_2\le R$ provided 
$|c(0)|=\|v_m(0)\|_2\le R$.  Thus the map
$$
S: \R^m \ni c(0)  \mapsto c(T)\in\R^M
$$ 
maps $\mathbb B_R^m$ into itself.  Since it is also continuous, we conclude by Brouwer's theorem that for each $m\ge1$ there exists $v_m(0)$ such that $v_m(0)=v_m(T)$. 
We may then extend $v_m(t)$ to the half-line $[0,\infty)$ to a periodic function of period $T$. 
Clearly, by \eqref{eq:2.7},
\begin{equation}\label{eq:2.9}
\|v_m(t)\|_2\le R+ 2\int_0^t {\rm e}^{2\tau/h^2}\|f(\tau)\|_2d\tau, \quad t\ge 0.
\end{equation}
Moreover, from \eqref{eq:2.5}, \eqref{eq:2.6} and   the time-periodicity of $v_m(t)$ we see that 
\begin{equation}\label{eq:2.10}
\Int0t\|\nabla v_m(\tau)\|_2^2 d\tau \le h^4
\Int0{\ell\,T}\|f(\tau)\|_2^2 d\tau, \quad t\in (0,\ell\,T), \ell>0\,.
\end{equation}
Furthermore, integrating both sides of \eqref{eq:2.4} between 0 and $t$ and using \eqref{eq:2.2}, \eqref{eq:2.3}, \eqref{eq:2.9} and \eqref{eq:2.10} we show that, 
for each fixed $r$, the sequence  of functions $\{(v_m(t),\psi_r)\}$ is  uniformly continuous and  uniformly bounded. Combining this information with \eqref{eq:2.4}, \eqref{eq:2.9}--\eqref{eq:2.10}, using \eqref{eq:2.2}--\eqref{eq:2.3}, 
and following the classical procedure (see, e.g., \cite[pp. 18--20]{Gal00}) we show the existence of a function $v\in L^\infty(0,\cT;\bH(\Omega))\cap L^2(0,\cT;\bH^1(\Omega))$ for all 
$\cT>0$  and a subsequence $\{v_{m'}\}$ such that for all $\cT>0$
\be\ba{lll}\label{eq:2.11}
%v\in L^\infty(0,\cT;\bH(\Omega)\cap L^2(0,\cT;\bH^1(\Omega))\\
v_{m'}&\to v\ &\mbox{weakly in $L^2(0,\cT;\bH^1(\Omega))\,,$}\\ 
v_{m'}(t)&\to v(t)\ &\mbox{weakly in $\bH(\Omega)$ for all $t\in [0,\cT]$\,,}\\
v_{m'}&\to v\ &\mbox{strongly in $L^2(0,\cT;\bH(\Omega))$}\,.
\ea
\ee
Recalling that $v_m(t+T)=v_m(T)$ for all $t\ge 0$, the second relation in \eqref{eq:2.11} implies that $v$ satisfies both properties i) and iv) of weak solutions. 
Furthermore, again by \eqref{eq:2.11} and \eqref{eq:2.5},  we see  that $v$ satisfies also property iii). Finally, we integrate \eqref{eq:2.4} over $(0,t)$ and then pass to 
the limit $(m') \to\infty$. Using \eqref{eq:2.11} and \eqref{eq:2.2}--\eqref{eq:2.3} we see that
$$
(v(t)-v(0),\psi_r)=\int_0^t(v\cdot\NP \psi_r,v)+(w(v)\,{\partial_z}\psi_r,v)+(\nabla v,\nabla\psi_r)+(f,\psi_r)\,,
$$
for all $r\ge 1$. From this equation, taking into account the mentioned properties of $\{\psi_n\}$, again by classical arguments (e.g., \cite[\S2]{Gal00}) we prove that 
$v$ satisfies also property ii), which concludes the proof. 
\QED

\section{Existence of Global Strong Solutions to the Initial-Value Problem: Inhomogeneous Case}

As mentioned earlier on, the basic idea for the proof of Theorem \ref{thm:main} is to compare the weak $T$-periodic solution of Proposition  \ref{prop:3.1} with a strong global solution $(u,p)$ to the inhomogeneous {\em initial-value problem}
\eqref{eqp.2}--\eqref{eqp.3}, for an appropriate choice of the initial data $a$. As customary \cite{GMR01,CT07}, by ``strong'' we mean  $u\in C([0,\cT]; \bH^1(\Omega)) \cap L^2(0,\cT; \bH^2(\Omega))$ with $\partial_t u\in L^1(0,\cT;\mathbb H(\Omega))$, 
and $p\in L^1(0,\cT; H^1(\Omega))$, $\cT>0$.  In \cite{HK15},  existence of such solutions has been established when $f\equiv 0$. In the following proposition, we shall suitably adapt the arguments of \cite{HK15} to prove existence of global strong solutions  to \eqref{eqp.2}--\eqref{eqp.3}, when $f\not\equiv 0$ is prescribed in a proper function class.  

\begin{proposition}\label{prop:4.1}
Let $\cT>0$ arbitrary, $a\in \bH^1(\Omega)$ and $f\in L^2(0;\cT; L^2(\Omega))$. Then, problem  \eqref{eqp.2}--\eqref{eqp.3} has a unique strong solutions in the interval $(0,\cT)$. 
\end{proposition}

\begin{proof}
The existence of a unique strong local solution was already proved in \cite[Theorem 1.2]{GMR01}. Hence, in order to prove the assertion, it suffices to show that the velocity field $u$, 
of a given local, strong  solution to the above problem, admits an a priori bound in the space 
$C([0,\cT]; \bH^1(\Omega)) \cap L^2(0,\cT; \bH^2(\Omega))$. Observe that this will also imply the stated properties on $\partial_tu$ and $p$, since from \eqref{eqp.2}$_1$ we first readily infer
$$
\|\partial_tu\|_2 \le \|u\|_\infty\|\nabla u\|_2+\|\nabla u\|_4^2+\|\Delta u\|_2\,,\ \ \|\nabla_Hp\|_2\le \|\partial_tu\|_2 +\|u\|_\infty\|\nabla u\|_2+\|\nabla u\|_4^2+\|\Delta u\|_2\,, 
$$
and then use the embedding $H^2(\Omega)\subset W^{1,4}(\Omega)\subset L^\infty(\Omega)$. 

Noticing that $\|u\|_{H^2} \le C\|\Delta u\|_2$,  in order to show the above bound, it suffices to  prove that
\begin{equation} \label{eq:LinfH1 cap L2H2 estimate}
	\|u(t)\|_{H^1}^2 + \int_0^t \|\Delta u(\tau)\|_{2}^2\,d\tau \le B(\|a\|_{\bH^1}, \|f\|_{L^2_{\cT}(L^2)}, \cT), \quad t\in [0,\cT],
\end{equation}
where $\|f\|_{L^2_{\cT}(L^2)} := (\int_0^{\cT} \|f(\tau)\|_2^2\,d\tau)^{1/2}$ and $B$ is a continuous function. 
Here and hereafter, $C$ denotes a generic constant.
	
In what follows, we shall closely employ the strategy of \cite{HK15} and show how the main estimates obtained there in (6.5), (6.7), (6.9) and (6.10) modify if $f\not\equiv 0$ satisfies the stated assumptions. This will be achieved in {\em Steps 1--4} below, which will then lead to the proof of \eqref{eq:LinfH1 cap L2H2 estimate} in {\em Step 5}. We begin to observe that multiplying both sides of \eqref{eqp.1} (written for $(u,p)$) by $u$, integrating by parts and using Schwartz inequality and Poincar\'e inequality \eqref{eq:2.6} we  show
\begin{equation} \label{eq:energy inequality}
	\|u(t)\|_{2}^2 + \int_0^t \|\nabla u(\tau)\|_{2}^2\,d\tau \le  \|a\|_{2}^2 + h^2\|f\|_{L^2_{\cT}(L^2)}^2,  \quad  t\in [0,\cT].
\end{equation}
The functions $\bar u = \frac1h \int_{-h}^0 u\,dz$ and $\tilde u := u - \bar u$ fulfill the following {equations:}
\begin{align}
	\partial_t \bar u - \Delta_H \bar u + \nabla_Hp &= \bar f -\bar u\cdot\nabla\bar u - \frac1h \int_{-h}^0 (\tilde u\cdot\nabla_H\tilde u + \mathrm{div}_Hu\,\tilde u)\,dz - 
\frac1h u_z|_{\Gamma_b} & \text{in}\quad G, \label{eq:vertical averaged PEs} \\
	\mathrm{div}_H\bar u &= 0 & \text{in}\quad G, \notag
\end{align}
with $u_z := {\partial_z} u$, as well as
\begin{equation} \label{eq:residual PEs}
	\partial_t\tilde u - \Delta\tilde u + \tilde u\cdot\nabla_H\tilde u + u_3u_z + \bar u\cdot\nabla_H\tilde u  = \tilde f - \tilde u\cdot\nabla_H\bar u + \frac1h \int_{-h}^0 (\tilde u\cdot\nabla_H\tilde u + \mathrm{div}_Hu\,\tilde u)\,dz + \frac1h u_z|_{\Gamma_b} \quad\text{in }\;\Omega.
\end{equation}

\noindent
\emph{Step 1:} Equation  \eqref{eq:vertical averaged PEs} implies
$$\begin{array}{rl}\medskip
	\partial_t\|\nabla_H\bar u(t)\|_{L^2(G)}^2 + &\|\Delta_H\bar u\|_{L^2(G)}^2 + \|\nabla_Hp\|_{L^2(G)}^2 
\\	&\le\; C_1(\big\| |\bar u||\nabla_H\bar u| \big\|_{L^2(G)}^2 + \big\| |\tilde u||\nabla_H\tilde u| \big\|_2^2 + \|u_z\|_{L^2(\Gamma_b)}^2 + \|f\|_2^2).
\end{array}
$$
Therefore, if $f\not\equiv 0$ estimate (6.5) in \cite{HK15} is replaced by 
\begin{align*}
\partial_t\|\nabla_H\bar u\|_{L^2(G)}^2 + \|\nabla_Hp\|_{L^2(G)}^2 &\le C(\|u\|_2 + \|u\|_2^2)(\|u\|_{H^1} + \|u\|_{H^1}^2)\|\nabla_H\bar u\|_{L^2(G)}^2 + 
C_1\big\| |\tilde u||\nabla\tilde u| \big\|_{L^2(G)}^2 \\
& \quad + \frac14 \|\nabla u_z\|_{L^2(G)}^2 + C(1+\|u\|_2^2 + \|u\|_2^4)\|u{\|_{H^1}}^2 + C\|f\|_2^2.
\end{align*}

%\vspace{.1cm}
\noindent
\emph{Step 2:} Multiplying \eqref{eqp.2} by $-{\partial_z} u_z$ and integrating by parts leads to
\begin{equation*}
	\frac12 \partial_t \|u_z\|_2^2 + \|\nabla u_z\|_{L_2(\Omega)}^2 
	= -\int_G \nabla_Hp\cdot u_z|_{\Gamma_b}  - \int_\Omega (u_z\cdot\nabla_Hu)\cdot u_z  + \int_\Omega \mathrm{div}_H u\, u_z\cdot u_z - \int_\Omega f\cdot {\partial_z}u_z.
\end{equation*}
Therefore, by using Cauchy-Schwartz inequality in the latter and proceeding as in \cite{HK15}, we obtain that equation (6.7) in \cite{HK15} generalizes to the following one
\begin{equation*}
	\partial_t\|u_z\|_2^2 + \|\nabla u_z\|_2^2 \le C\, \big(\|u\|_{H^1} + \|u\|_{H^1}^2\big)\|u_z\|_2^2 + 
\frac12\|\nabla_Hp\|_{L^2(G)}^2 + 2C_2\big\| |\tilde u||\nabla\tilde u| \big\|_2^2 +C\big( \|u\|_{H^1}^2 + \|f\|_2^2\big).
\end{equation*}

%\vspace{.1cm}
\noindent
\emph{Step 3:}
Equation (\ref{eq:residual PEs}) implies
\begin{align*}
	\frac14\partial_t\|\tilde u\|_4^4 + \frac12\big\|\nabla |\tilde u|^2\big\|_2^2 + \big\| |\tilde u||\nabla\tilde u| \big\|_2^2  &= - 
\int_\Omega (\tilde u\cdot\nabla_H\bar u) \cdot |\tilde u|^2\,\tilde u 
	 + \frac1h \int_\Omega \int_{-h}^0 (\tilde u\cdot\nabla_H\tilde u 
+ \mathrm{div}_Hu\,\tilde u)\,dz \cdot |\tilde u|^2\,\tilde u \\ & \quad + 
\frac1h \int_\Omega u_z|_{\Gamma_b}\cdot |\tilde u|^2\,\tilde u + \int_\Omega \tilde f\cdot |\tilde u|^2\,\tilde u.
\end{align*}
The last term on the right-hand side is bounded by $\|\tilde f\|_2 \big\| |\tilde u|^3 \big\|_2$. We see moreover that
\begin{align*}
	\big\| |\tilde u|^3 \big\|_2 &= \big\| |\tilde u|^2 \big\|_3^{3/2} \le C \big\| |\tilde u|^2 \big\|_{W^{1/2,2}(\Omega)}^{3/2} 
		\le C \big\| |\tilde u|^2 \big\|_2^{3/4} \Big(\big\| |\tilde u|^2 \big\|_2 + \big\| \nabla|\tilde u|^2 \big\|_2 \Big)^{3/4} \\
&= C \| \tilde u \|_4^{3/2} \Big(\big\| |\tilde u|^2 \big\|_2 + \big\| \nabla|\tilde u|^2 \big\|_2 \Big)^{3/4} 
		\le C \|\tilde u\|_4^3 + C \|\tilde u\|_4^{3/2} {\big\| \nabla|\tilde u|^2 \big\|_2^{3/4}},
\end{align*}
where we have used the embedding $W^{1/2,2}(\Omega) \hookrightarrow L^3(\Omega)$ and an interpolation inequality.
It then follows that
\begin{align*}
	\int_\Omega \tilde f\cdot |\tilde u|^2\,\tilde u &\le C\|\tilde f\|_2 \|\tilde u\|_4^3 + C\|\tilde f\|_2 \|\tilde u\|_4^{3/2} \big\| \nabla|\tilde u|^2 \big\|_2^{3/4} \\
		&= C \big(\|\tilde f\|_2^2 \|\tilde u\|_4^4\big)^{1/2} \big(\|\tilde u\|_4^2\big)^{1/2}  +  C \|\tilde f\|_2^{2/8}  \big(\|\tilde f\|_2^2 \|\tilde u\|_4^4\big )^{3/8}  \big( \big\| \nabla|\tilde u|^2 \big\|_2^2 \big)^{3/8} \\
		&\le \frac1{{6}}\big\| \nabla|\tilde u|^2 \big\|_2^2 + C\|\tilde f\|_2^2 \|\tilde u\|_4^4 + C\|\tilde u\|_4^2 + C\|\tilde f\|_2 \\
		&\le \frac1{{6}}\big\| \nabla|\tilde u|^2 \big\|_2^2 + C\|f\|_2^2 \|\tilde u\|_4^4 + C\|u{\|_{H^1}}^2 + C\|f\|_2.
\end{align*}
Therefore, estimate (6.9) in \cite{HK15} is now replaced by 
\begin{align*}
	\frac{C_3}{4} \partial_t\|\tilde u\|_4^4 + C_3 \big\||\tilde u||\nabla\tilde u|\big\|_2^2 
	\le\;  C\big( \|u\|_{H^1}^{2/3} + \|u\|_{H^1} + \|u\|_{H^1}^2 + \|f\|_2^2 \big) \|\tilde u\|_4^4 + \frac14 \|\nabla u_z\|_2^2 + C\|u\|_{H^1}^2 + C\|f\|_2,
\end{align*}
where $C_3 := 2(C_1 + 2C_2)$.

\vspace{.2cm}\noindent
\emph{Step 4:} Combining {\em Steps 1--3}, we see that estimate (6.10) in \cite{HK15} is now being replaced by 
\begin{align*}
	\partial_t \big( 8 \|\nabla_H\bar u\|_{L^2(G)}^2 & + \|u_z\|_2^2 + \frac{C_3}{4}\|\tilde u\|_4^4\big) + 
\frac12 \big(\|\nabla_Hp\|_{L^2(G)}^2 + \|\nabla u_z\|_2^2 + C_3\big\| |\tilde u||\nabla_H\tilde u| \big\|_2^2\big) \\
&	\le K_1(t)\big(8\|\nabla_H\bar u\|_{L^2(G)}^2 + \|u_z\|_2^2 + (C_3/4)\|\tilde u\|_4^4\big) + K_2(t),
\end{align*}
where, thanks to (\ref{eq:energy inequality}), the functions $K_1$ and $K_2$ given by  
\begin{align*}
	K_1(t) &:= C(1+\|u\|_2+\|u\|_2^2) (\|u\|_{H^1}^{2/3} + \|u\|_{H^1} + \|u\|_{H^1}^2 + \|f\|_2^2), \\
	K_2(t) &:= C(1+\|u\|_2^2 + \|u\|_2^4)\|u\|_{H^1}^2 + C(\|f\|_2 + \|f\|_2^2),
\end{align*}
are integrable on $[0, \cT]$.
We thus conclude by  Gronwall's inequality that
\begin{equation} \label{eq:key a priori estimate}
	\|\nabla_H\bar u(t)\|_{L^2(G)} + \|u_z(t)\|_2 + \|\tilde u(t)\|_4 + \int_0^t \Big(\|\nabla_Hp\|_{L^2(G)}^2 + \|\nabla u_z\|_2^2 + \big\| |\tilde u||\nabla\tilde u| \big\|_2^2 \Big) \,d\tau
\end{equation}
must be bounded by some $B_1(\|a{\|_{H^1}}, \|f\|_{L^2_{\cT}(L^2)}, \cT)$ for all $t\in[0,\cT]$.

\vspace{.2cm}\noindent
\emph{Step 5:}
Following the estimates in Step 5 of  \cite{HK15}, we obtain now 
\begin{equation*}
	\partial_t\|\nabla u\|_2^2 + \|\Delta u\|_2^2 \le C\|u_z\|_2^2\|\nabla u_z\|_2^2\|\nabla u\|_2^2 + C(\|\bar u\|_{H^1(G)}^2 + \|\tilde u\|_4^4)\|u\|_{H^1}^2 + C\|f\|_2^2,
\end{equation*}
where, thanks to (\ref{eq:energy inequality}) and to the estimate for (\ref{eq:key a priori estimate}),
\begin{align*}
	C\|u_z\|_2^2\|\nabla u_z\|_2^2 \quad\text{and}\quad C(\|\bar u\|_{H^1(G)}^2 + \|\tilde u\|_4^4)\|u\|_{H^1}^2 + C\|f\|_2^2
\end{align*}
are integrable on $[0, \cT]$.
Thus, Gronwall's inequality yields the desired estimate (\ref{eq:LinfH1 cap L2H2 estimate}). The proof is complete.
\end{proof}

\section{Weak-Strong Uniqueness and Proof of the Main Result}
In this final section we give a proof of Theorem \ref{thm:main}  based on a weak-strong uniqueness argument for the {\em initial-value} problem.

More precisely, given $f \in L^2(0,T;L^2(\Omega))$, let $v$ be a weak $T$-periodic solution corresponding to Proposition \ref{prop:3.1}. 
By properties i) and iii) we infer that there exists $t_0>0$ such that $v(t_0)\in \bH^1(\Omega)$ and   
\be\label{eq:4.1}
\|v(t)\|_2^2 + 2\int_{t_0}^t\|\nabla v(\tau)\|_2^2d\tau\le \|v(t_0)\|_2^2+2\int_{t_0}^t(f(\tau),v(\tau))d\tau, \quad t\ge t_0, 
\ee
while from ii)  we deduce that, for arbitrary $\cT>0$, the following relation 
\be\label{eq:4.2}
\Int{t_0}t\big\{(v,\partial_t\varphi) - (\nabla v,\nabla\varphi)-(v\cdot\NP v+w(v)\,{\partial_z} v,\varphi)\big\}
=-\Int{t_0}t(f,\varphi)+(v(t),\varphi(t))-(v(t_0),\varphi(t_0))
\ee
holds for all $t\in (t_0,\cT]$ and all $\varphi\in C^1([t_0,\cT];\bH^1(\Omega))\cap L^2(0,\cT;\bH^2(\Omega))$. 
Note that in  \eqref{eq:4.2} we have used the  identity
\be\label{eq:4.3}
({\sf v}\cdot\NP \phi+w({\sf v})\,{\partial_z} \phi,{\sf v})=-({\sf v}\cdot\NP{\sf v}+w({\sf v})\,{\partial_z} {\sf v},\phi), \quad ({\sf v},\phi)\in \bH^1(\Omega)\times \bH^2(\Omega), 
\ee
which follows by integration by parts.

We now look at our weak time-periodic solution as  a weak solution to the {\em initial-value problem} with initial data $a\equiv v(t_0)$. Since $v(t_0)\in \bH^1(\Omega)$ 
we may use it also as { initial value} for the {\em global} strong solution determined in Proposition \ref{prop:4.1}. The assertion of Theorem \ref{thm:main}  follows provided  
we are able to show that the weak  solution  coincides with the strong one. 

As a matter of fact, such a weak-strong uniqueness result is already known, under slightly different boundary conditions; see \cite[Theorem 1.3]{GMR01}.  The arguments used there 
would equally apply to the case at hand. Nevertheless, we shall sketch a proof here. To this end, let $u$ be the velocity field of the strong solution determined in  Proposition \ref{prop:4.1} 
corresponding to $f$ and initial data $v(t_0)$. Clearly, $u$ satisfies
\be\label{eq:4.4}
\Int{t_0}t\big\{(u,\partial_t\varphi) -(\nabla u,\nabla\varphi) {+} (u\cdot\NP \varphi+w(u)\,{\partial_z} \varphi,u)\big\}
=-\Int{t_0}t(f,\varphi)+(u(t),\varphi(t))-(v(t_0),\varphi(t_0)), 
\ee
for all $t\in (t_0,\cT]$ and  all  $\varphi$.
In addition, thanks to the regularity properties of $u$, we see that
\be\label{eq:4.5}
\|u(t)\|_2^2 + 2\int_{t_0}^t\|\nabla u(\tau)\|_2^2d\tau= \|v(t_0)\|_2^2+2\int_{t_0}^t(f(\tau),u(\tau))d\tau, \quad t\ge t_0. 
\ee
Next, let
$$
v_h(t):=\int_0^\cT j_h(t-s)v(s)\,ds \mbox{ and } u_h(t):=\int_0^\cT j_h(t-s)u(s)\,ds
$$
be the (Friedrichs) time-mollifier of $v$ and $u$, respectively, where $j_h\in C_c^\infty(-h,h)$, $0<h<\cT$, is even and positive with $\int_\R j_h(s)ds=1$. Then, as is well known,
%\begin{align}
\be
\ba{lll}
\lim_{h\to 0} \int_0^\cT\|v_h(\tau)-v(\tau)\|^2_{\bH^1}&=0, \qquad \essup{t\in[0,\cT]}\,\|v_h(t)\|_2 &\le \essup{t\in[0,\cT]}\,\|v(t)\|_2, \mbox{ and} \\ % \label{eq:4.6}\\
\lim_{h\to 0} \int_0^\cT\|u_h(\tau)-u(\tau)\|^2_{H^2}&=0, \qquad \essup{t\in[0,\cT]}\,\|u_h(t)\|_{H^1} &\le \essup{t\in[0,\cT]}\,\|u (t)\|_{H^1}. \label{eq:4.7}
\ea
\ee
%\end{align}
Moreover, 
\be\label{eq:4.8}
\Int{t_0}t(v,\partial_tu_h)=-\Int{t_0}t(u,\partial_tv_h)\,,
\ee
while the weak continuity of $v$ and $u$ implies 
\be\label{eq:4.9} 
\lim_{h\to 0}(u(t),v_h(t))=\lim_{h\to 0}(u_h(t),v(t))=\half(u(t),v(t))\,,\ \ t\ge t_0.
\ee
Finally, setting $\sigma:=v-u$, by integrating by parts we get 
\be\label{eq:4.10}
(\sigma\cdot\NP u+w(\sigma){\partial_z} u,u)=0\,,\ \ \sigma\in {\bH^1}(\Omega)\,.
\ee
We now replace $\varphi$ in \eqref{eq:4.2} by $u_h$,  then use \eqref{eq:4.3} with ${\sf v}\equiv u$ and $\phi\equiv \varphi$ in \eqref{eq:4.4} and replace in the latter $\varphi$ by $v_h$.  Summing side by side the resulting equations and employing \eqref{eq:2.2}, \eqref{eq:2.3},  
\eqref{eq:4.7}, \eqref{eq:4.9} and \eqref{eq:4.10} we show that
\be\label{eq:4.11}
\Int{t_0}{t}\big\{-2(\nabla v,\nabla u) -(\sigma\cdot\NP\sigma+w(\sigma){\partial_z}\sigma,u)\big\}\,d\tau =-\Int{t_0}t f\cdot(u+v)\,d\tau-\big[(v(t),u(t))-\|v(t_0)\|_2^2\big]\,,
\ee
Adding twice \eqref{eq:4.11} to  \eqref{eq:4.1} and \eqref{eq:4.5} we infer
$$
\|\sigma(t)\|_2^2+2\int_0^t\|\nabla\sigma(\tau)\|_2^2\,d\tau\le 2\int_{t_0}^t\{(\sigma(\tau)\cdot\NP\sigma(\tau)+w(\sigma){\partial_z}\sigma(\tau),u(\tau))\big\}d\tau\,.
$$
In order to estimate the above right-hand side we use \eqref{eq:4.3} and \eqref{eq:2.2}, \eqref{eq:2.3} with $v_1=v_2\equiv\sigma$ and $v_3\equiv u$, respectively, and Poincar\'e's and Young's inequalities to deduce
\begin{equation*}
\ba{rl}%\label{eq:4.12}
\|\sigma(t)\|_2^2 + 2 \Int{t_0}t\|\nabla\sigma(\tau)\|_2^2\,d\tau &\le C\!\Int{t_0}t\{\|\sigma(\tau)\|_2^{\frac12} \|\sigma(\tau)\|_{H^1}^{\frac32} (\|\NP u(\tau)\|_2
+ \|\partial_z u(\tau)\|_2^{\frac12} \|\partial_z u(\tau)\|_{H^1}^{\frac12} )\}d\tau\\
&\le C \Int{t_0}t g(\tau) \|\sigma(\tau)\|_2^2 d\tau + \Int{t_0}t\|\nabla\sigma(\tau)\|_2^2 d\tau\,, 
\ea
\end{equation*}
where $g(\tau):=\|\NP u(\tau)\|^4_2 + \|\partial_z u(\tau)\|_2^2 \|\partial_z u(\tau)\|_{H^1}^2$. By the regularity properties of $u$ it follows that $g\in L^1(t_0,t)$ for  arbitrary $t> t_0$, 
so that by Gronwall's lemma we conclude that 
$$
\sigma\equiv v-u=0 \mbox{ a.e. in } \Omega\times (t_0,t) \mbox{ for all } t>t_0.
$$ 
By the time-periodicity of $v$  this proves the desired regularity for $v$. The proof of the Theorem \ref{thm:main} is complete. 
%\QED

\end{document}